\documentclass[12pt,a4paper]{amsart}
\usepackage{amsfonts,amssymb,amscd,amsmath,enumerate,verbatim,calc}
\usepackage{mathrsfs}
\headsep=1cm \numberwithin{equation}{section}
\newtheorem{thm}{Theorem}[section]
\newtheorem{cor}[thm]{Corollary}
\newtheorem{lem}[thm]{Lemma}
\newtheorem{ex}[thm]{Example}
\newtheorem{prop}[thm]{Proposition}
\theoremstyle{definition}

\theoremstyle{Remark}

\newtheorem{rem}[thm]{Remark}

\numberwithin{equation}{section}
\usepackage[linktoc=all, hyperindex]{hyperref}% 
\hypersetup{colorlinks, citecolor=blue, filecolor=blue, linkcolor=red, urlcolor=black}
\newcommand\diam{\operatorname{diam}}

\begin{document}
\title[On Commuting graphs of triangular rings]{On Commuting graphs of triangular rings }
\author[]{H. CHERAGHPOUR$ ^{*} $, M.N. Ghosseiri, M. Jafari AND  F. SEYFPOUR}
%\address{This work is supported in part by the Slovenian Research Agency, project N1-0210.} 
\address{University of Primorska, FAMNIT and IAM, Glagolja{\v s}ka 8, 6000 Koper, Slovenia.}
\email{cheraghpour.hassan@yahoo.com}

\address{Department of Mathematics, University of Kurdistan, P.O. Box 416, Sanandaj, Iran.}
\email{mnghosseiri@yahoo.com}
\email{madineh.jafari3978@gmail.com}
\email{farnazs1364@gmail.com}

%\subjclass{}%
\thanks{This work is supported in part by the Slovenian Research Agency (project N1-0210).}
\thanks{$ ^{*} $Corresponding Author: Hassan Cheraghpour.}

%\date{}%
%\dedicatory{}%
%\commby{}%
% ----------------------------------------------------------------
\maketitle
% 1-----------------------------------------------------------------------------------------------

\begin{abstract} 
Let $ R $ be a noncommutative ring with identity. The commuting graph of $ R $, denoted by $ \Gamma(R) $, is a graph with vertex set $ R \setminus Z(R) $, and two vertices $ a $, $ b $ are adjacent if $ a\neq b $ and $ ab=ba $. Let $ T=Tr(R) $ be the ring of all $ 2\times 2 $ upper triangular matrices over $ R $ and $ \Gamma(T) $ be the commuting graph of $ T $. In this article, we find the number of edges, cliques, clique number, and independence number of $ \Gamma (T) $ when $ R $ is a finite field. Moreover, we show that for the case when $ R= \mathbb{Z}_{n} $ is not a field, $ \Gamma (T) $ is connected with diameter 3. Some useful related results are also obtained, some examples are presented and a question is posed.
\\ 
\end{abstract}
\maketitle 
% 1-----------------------------------------------------------------------------------------------

\section{introduction}
Let $ R $ be a noncommutative ring with center $ Z(R) $. The commuting graph of $ R $ denoted by  $ \Gamma(R) $, is a graph with vertex set $ R \setminus Z(R) $ and two distinct vertices $ a $, $ b$ are adjacent if $ ab=ba $.
First, this graph was introduced by Akbari et al. \cite{AGHM} in 2004, for classification of rings using their commuting graphs.
Then, in 2006, Akbari and Raja \cite{AR} showed that in some conditions, the graphs associated with upper triangular matrices and full matrix rings are connected.
 Some times later, Akbari et al. \cite{AMRR} investigated the diameters of $ \Gamma (M_{n}(D)) $ and determined the diameters of some induced subgraphs of $\Gamma( M_{n}(D)) $, where $ D $ is a divison ring.

Next, in 2008, Abdollahi \cite{AA} showed that for any noncommutative ring $ R $ and finite field $ F$, if $ \Gamma(R) \cong \Gamma (M_{n}(F)) $, then $ |R|=|M_{n}(F)| $. Then, some times later, Akbari et al. \cite{ABM} investigated some graph-theoretic properties of $ \Gamma (M_{n} (F)) $, where $ F $ is a field and $ n\geqslant 2 $.

Next, in 2010, Akbari et al. \cite{AKR} investigated some graph-theoretic properties of $ \Gamma (KG) $, where $ G $ is a finite group, $ K $ is a field and $ 0\neq |G| \in K $. Then, in 2010, Mohammadian \cite{AM} investigated the commuting graphs of finite matrix rings. Next, Giudici and Pope in \cite{GP} studied the connectivity and diameters of matrix rings over $\mathbb{Z}_{m} $.
In 2011, Omidi and Vatandoost \cite{OV} characterized all rings where the complements of their commuting graphs are planar. Then, in 2012, Dolinar et al. \cite{DKO} characterized the diameters of commuting graph of $ \Gamma (M_{n}(F)) $, where $ n\geq 3 $ and $ F $ is a field.
Next, in 2013 Ambrozie et al. \cite{ABKM} studied the connectedness of $ \Gamma (B(H)) $ in the case when $ H $ is an infinite-dimensional complex Hilbert space. Then, in 2014 Dolinar et al. \cite{DGKO} determined the conditions for matrix centeralizers which can guarantee the connectedness of the commuting graph for the full matrix algebra $ M_{n}(F) $ over an arbitrary field $ F $.
Then, in 2016, Dolzan et al. \cite{DKKO} partially proved the conjecture stated by Akbari and Raja \cite{AR} in 2006.
Some times later, Giudici and Kuzma \cite{GK} obtained complete classifications of the graphs with an isolated vertex or edge that are the commuting graph of a group and the cycles that are commuting graph of a centerfree semigroup.
Then, Vatandoost and Ramezani \cite{VR} investigated noncommutative rings of order $ p^{n} $, where $ p $ is prime and $ n\in \lbrace 4,5\rbrace $.
Finally in 2018, Dolzan et al. \cite{DKK1} showed that the diameter of the commuting graph of $ M_{n}(F) $, where $ n=p^{2} $ is an odd prime square and $ F $ is a finite field with sufficiently large number of elements, is at least five.\\
\indent Let $ G=(V,E) $ be a graph, where $ V=V(G) $ is the set of all vertices of $ G $ and $ E=E(G) $ is a set of unordered pairs of vertices (edges) of $ G $.
For each $ v\in V $, define the neighbourhood of $ v $ by
\begin{align}\label{}
N(v)= \lbrace u \in V | (u, v) \in E \rbrace,
\end{align}
and define the {\it degree} of $ v $ by $ d(v)= |N(v)| $. The minimum and the maximum degrees of the vertices of $ G $ are denoted by $ \delta (G) $ and $ \Delta (G) $, respectively. $ G $ is $ d ${\it-regular} if all vertices have degree $ d $. A  walk in $ G $ is a sequence of vertices $v_{0}, v_{1}, \cdots , v_{k} $ and a sequence of edges $ (v_{i},v_{i+1}) \in E(G) $. A  walk is called a {\it path} if all $ v_{i} $ are distinct. If for such a  path (with $ k\geqslant 2$), $ (v_{0},v_{k}) $ is also an edge in $ G $, then $v_{0},v_{1}, \cdots ,v_{k} $ is called a  cycle. The {\it length} of a path,  cycle or walk, is the number of edges in it. The graph $ G $ is {\it connected} if for all pairs $ u,v \in V(G) $, there is a path from $ u $ to $ v $. (Note that, by definition, it suffices for there to be a  walk from $ u $ to $ v $.)
A graph $ (H,F) $ is a {\it subgraph} of $ (G,E) $ if $ H\subseteq G $ and $ F\subseteq E $. The subset $ X \subseteq V $ is called a {\it clique} if the induced  subgraph by $ X $ (whose edges are all $ (a, b) $ in $ E $ such that $ a,b \in X $) is a complete graph. The clique number $ \omega (G) $ of $ G $ is the maximum size of a clique in $ G $.
The subset $ X \subseteq V $ is called an {\it independent set} if the induced subgraph on $ X $ has no edges; i.e, no pairs in $ X $ make an edge. The maximum size of such an $ X $ is called the {\it independence number} of $ G $ and denoted by $ \alpha (G) $. A (connected) {\it component} of $ G $ is a connected subgraph of $ G $ which is maximal with respect to inclusion. The complement $\overline{G} $ of $ G=(V,E) $ is the graph with the same vertex set $ V $, but $ (u,v) \in E (\overline{G} ) $ if and only if $ (u,v) \notin E=E(G)$. Let $ u, v \in V $ and let $ d(u,v) $ be the length of shortest  path between $ u $ and $ v $, if such a path exists. Otherwise, define $ d(u,v)=\infty $. The {\it diameter} of $ G $ is defined by 
\[
\diam G=\sup \lbrace d(u,tv)\, |\, u, v \text{ are distinct vertices of}\,\, G \rbrace.
\]

The paper is organized as follows. After this introductory section, we state the main results and present some examples in Section 2.\\

Let $ R $ be a ring with identity, $ T=Tr (R) $ be the ring of all $ 2\times 2 $ upper triangular matrices over $ R $, and let $ \Gamma (T) $ be the commuting graph of $T $. In this paper we find the number of edges, cliques, clique number, independence number, and the {\it diameter} of $ \Gamma (T) $ for the case when $R$ is a finite feild (Theorem \ref{thm2.9}). Moreover, we show that for the case when $T= \mathbb{Z}_{n}$ is not a field, $ \Gamma (T) $ is connected, $ \diam \Gamma (T)=3 $ (Theorem \ref{thm2.10}), and pose a question on the independence number of $ \Gamma (T) $ (see page 8). Further, for the case when $ R= \mathbb{Z}_{p^{n}} $, where $ p $ is a prime and $ n\geqslant 2 $, we find the number of edges, clique number, and the independence number of $ \Gamma (T) $ (Corollary \ref{cor2.3}). \\

We begin by the following obvious, but crucial remark which will be used frequently in the sequel. \\

\begin{rem} \label{1.1}
Let $ R $ be an arbitrary ring with identity and center $ Z(R) $. For any $ x \in R $, let $ C_{R} (x) $ be the centralizer of $ x $ in R. Then, for every $ x \in R $ and $ z \in Z(R) $ we have 
\begin{align} \label{2} 
C_{R} (x)= C_{R} (x+z).
\end{align}
\end{rem}

Now, let $ T=Tr(R) $. Noting that $ Z(T)={ \lbrace aI | a \in Z(R)} \rbrace $, where $ I $ is the identity matrix, it follows from (\ref{2}) that, for any $ A = \begin{pmatrix} 
x & y\\
0 & z
\end{pmatrix} \in T $
and $ a \in Z(R) $, we have
\begin{align*}
C_{T} (A) & =C_{T}
\begin{pmatrix}
x-a & y\\
0 & z-a
\end{pmatrix}.
\end{align*} 
Therefore, in order to find the edges, neighbourhoods, degrees,..., of the commuting graph of $ Tr(R) $, where $ R $ is commutative, it suffices to work with the matrices of the form 
$ \begin{pmatrix}
x & y \\
0 & 0
\end{pmatrix} $
in $ Tr(R) $.\\

Assume that $ R $ is a finite commutative ring with identity and let $ T=Tr(R) $. Then the number of vertices of $ \Gamma (T) $ is 
\begin{align} \label{3}
|V(\Gamma (T)|=|T\setminus Z(T)| = |R|^{3} - |R|.
\end{align}

Recall that for any vertex $ A $ in $ \Gamma (T) $, the  neighbourhood of $ A $ is the set 
\begin{align} \label{4}
N(A)= \lbrace B\in T| AB=BA , A\neq B , B \notin Z(T) \rbrace ,
\end{align}
and the number of (unordered) edges ($ A,B $ ) in $ \Gamma (T) $ is the degree of $ A $= $ d(A)=|N(A)| $.
\\
In fact,
\begin{align} \label{5}
d(A)=\mid \lbrace B\in C_{T}(A) | B \neq A , B \neq aI ~for~ all ~ a \in R \rbrace \mid .
\end{align}

By the remark above, in order to find the number of edges $ | E(\Gamma (T))| $ of $ \Gamma (T) $, it suffices to consider the matrices of the following forms: 
\begin{align*}
A_{1} =
\begin{pmatrix}
u & 0 \\
0 & 0
\end{pmatrix},
A_{2} =
\begin{pmatrix}
0 & u \\
0 & 0
\end{pmatrix},
A_{3} =
\begin{pmatrix}
u & v \\
0 & 0
\end{pmatrix},
A_{4} =
\begin{pmatrix}
n & 0 \\
0 & 0
\end{pmatrix},\\
\end{align*}
\begin{align}\label{6}
\,\,\,\,A_{5} =
\begin{pmatrix}
0 & n \\
0 & 0
\end{pmatrix},
A_{6} =
\begin{pmatrix}
n_{1} & n_{2} \\
0 & 0
\end{pmatrix},
A_{7} =
\begin{pmatrix}
u & n \\
0 & 0
\end{pmatrix},
A_{8} =
\begin{pmatrix}
n & u \\
0 & 0
\end{pmatrix}, 
\end{align}
\noindent where $ u,v \in U( R ) $, and $ n $, $n_{1} $, $n_{2} $ are in the set 
\begin{align} \label{7}
X= \lbrace r\in R\, |\, r\neq 0\, \text{and}\, r \, \text{is a zero divisor} \rbrace .
\end{align}
\indent Note that for the case when $ R $ is a finite field of order $p^{n} $, where $ p $ is a prime and $ n\geqslant 1 $, we have $ X=\emptyset $ and if $ R= \mathbb{Z} _{n}, n\geqslant 2 $, we have $ |X |=n-\varphi (n)-1 $, where $ \varphi $ is the Euiler function.

\indent For each $ 1 \leqslant i \leqslant 8 $, define $ A_{i} ^{*} $ as follows:
\begin{align*}
A_{1} ^{*} =&\,\lbrace \begin{pmatrix}
u+z & 0\\
0 & z
\end{pmatrix} | u \in U( R ) , z \in R \rbrace, \\
A_{2} ^{*} =&\,\lbrace \begin{pmatrix}
z & u\\
0 & z
\end{pmatrix} | u \in U( R ) , z \in R \rbrace, \\
A_{3} ^{*} =&\,\lbrace \begin{pmatrix}
u+z & v\\
0 & z
\end{pmatrix} | u,v \in U( R ) , z \in R \rbrace, \\
A_{4} ^{*} =&\,\lbrace \begin{pmatrix}
n+z & 0\\
0 & z
\end{pmatrix} | n \in X , z \in R \rbrace, \\
A_{5} ^{*} =&\,\lbrace \begin{pmatrix}
z & n\\
0 & z
\end{pmatrix} | n \in X , z \in R \rbrace, \\
A_{6} ^{*} =&\,\lbrace \begin{pmatrix}
n_{1}+z & n_{2} \\
0 & z
\end{pmatrix} | n_i \in X , z \in R \rbrace, \\
A_{7} ^{*} =&\,\lbrace \begin{pmatrix}
u+z & n\\
0 & z
\end{pmatrix} |u \in U(R), n \in X , z \in R \rbrace, 
\end{align*}
\begin{align}\label{8}
\,A_{8} ^{*} =&\,\lbrace \begin{pmatrix}
n+z & u\\
0 & z
\end{pmatrix} | n \in X ,u \in U(R), z \in R \rbrace. 
\end{align}
\indent Note that for the case when $ R $ is a field, $ A_{i} ^{*}=\emptyset $ for all $i\geqslant 4 $. Moreover, it is not hard to observe that $ A_{i} ^{*^{,}}$s make a partition of $ T\setminus Z(T) $:

\begin{align}\label{9}
T\backslash Z(T) = \bigcup ^{8} _{i=1} A_{i} ^{*}.
\end{align} 
\indent The following lemma will be frequently used in the sequel. (See \cite{HV}.)
\begin{lem}[Handshaking Lamma]\label{10}
Let $ G=(V,E) $ be a finite graph. Then 
$$ \sum_{v\in V} d(v)=2|E|. $$
\end{lem}
%-----------------------------------------------------------------------------------------------
\section{Main results}
In this section we state the main results of the paper, present some examples to illustrate them, and pose a question.\\

Let $ A_i $ and $ A_i^{*} $ be as in (\ref{6}), (\ref{8}). Using  Remark \ref{1.1} and the identity (\ref{9}), we have 
\begin{lem}\label{lem2.1}
Let $ A=\begin{pmatrix}
x & y\\
0 & 0
\end{pmatrix} $ be an arbitrary element in $ A^{*}  _{i} $ for some $ 1\leqslant i \leqslant 8 $. Define $ d(A^{*} _{i} ):= \sum _{B \in A^{*} _{i}} d(B) $. Then
\begin{align}\label{(2.1)}
d(A^{*}_{i} )= |A^{*}_{i} | d(A).
\end{align}
\end{lem}
\proof
It suffices to note that for each $ 1\leqslant i \leqslant 8 $, $ N(B)=N(A) $ for all $ B \in A_{i} ^{*} $.
\qed\\

Lemmas \ref{10}, \ref{lem2.1}, and the following theorem play a key role for several subsequent computations.
\begin{thm}\label{thm2.2}
Let $ R $ be a finite commutative ring and let $ A _{i},A^{*} _{i} $ be as in (\ref{6}), (\ref{8}), where $ u, v \in U(R) $ and $ n, n_1, n_2 \in X $ are all fixed. For each $ r \in R $, let $ Z(r) = \lbrace s \in R ~| ~sr=0 \rbrace $. Then
\begin{enumerate} 
\item[$(i)$] $ N(A_1)= \lbrace B=\begin{pmatrix}
a & 0\\
0 & c
\end{pmatrix}~|~ a, c \in R, ~~ B \neq aI, ~~ B\neq A_1 \rbrace ; $\\
$ d(A_1)=|R|^2 -|R| -1. $

\item[$(ii)$] $ N(A_2)= \lbrace B=\begin{pmatrix}
a & b\\
0 & a
\end{pmatrix}~|~ a,b \in R, ~~ B \neq aI, ~~ B\neq A_2 \rbrace ; $\\
$ d(A_2)=|R|^2 -|R| -1=d(A_{1}). $

\item[$(iii)$] $ N(A_3)= \lbrace B=\begin{pmatrix}
a & b\\
0 & a-u^{-1}vb
\end{pmatrix}~|~ a,b \in R, ~~ B \neq aI, ~~ B\neq A_3 \rbrace ;  $\\
$ d(A_3)=|R|^2 -|R| -1. $

\item[$(iv)$] $ N(A_4)= \lbrace B=\begin{pmatrix}
a & b\\
0 & c
\end{pmatrix}~|~ a,c \in R,~~ b \in Z(n), ~~ B \neq aI, ~~ B\neq A_4 \rbrace ; $\\
$ d(A_4)=|Z(n)|~|R|^2 -|R| -1. $

\item[$(v)$] $ N(A_5)= \lbrace B=\begin{pmatrix}
a & b\\
0 & a+z
\end{pmatrix}~|~ a,b \in R,~~ z \in Z(n), ~~ B \neq aI, ~~ B\neq A_5 \rbrace ; $\\
$ d(A_5)=|Z(n)|~|R|^2 -|R| -1 =d(A_{4}). $

\item[$(vi)$]
$ N(A_6)= \lbrace  B=\begin{pmatrix}
a & b\\
0 & c
\end{pmatrix}~|~ a,b,c \in R,~\text{such~that}\\
\indent  \indent ~~~~~~~~\,\,\,\,\,\,\,\,\,\,\,\,\,\,\, bn_1+(c-a)n_2=0, ~~ B \neq aI, ~~ B\neq A_6 \rbrace ,$\\
$ d(A_6)=S~|R|^2 -|R| -1, $ where $ S $ is the number of all $ x,y \in R $ that satisfy the equation $ n_1x+n_2 y =0 $.
\item[$(vii)$] $ N(A_7)= \lbrace B=\begin{pmatrix}
a & u^{-1}(a-c)n\\
0 & c
\end{pmatrix}~|~ a,c \in R, ~~ B \neq aI, ~~ B\neq A_7 \rbrace ; $\\
$ d(A_7)=|R|^2 -|R| -1. $

\item[$(viii)$] $ N(A_8)= \lbrace B=\begin{pmatrix}
a & b\\
0 & a-u^{-1}bn
\end{pmatrix}~|~ a,b \in R, ~~ B \neq aI, ~~ B\neq A_8 \rbrace ; $\\
$ d(A_8)=|R|^2 -|R| -1 =d(A_{7}). $
\end{enumerate}
\end{thm}
\proof
($ i $) Let 
$A_{1} = \begin{pmatrix}
u & 0\\
0 & 0
\end{pmatrix}$ $ \in A_{1} ^{*} $,  and assume that 
$ B=\begin{pmatrix}
a & b\\
0 & c
\end{pmatrix}$
$ \in T=Tr(R) $
satisfies $ AB=BA, A\neq B $, and 
$ B \notin \lbrace aI |a \in R \rbrace=Z(T) $. Since $ R $ is commutative, from $ AB=BA $ we conclude that $ ub=0 $, so that $ b=0 $. Therefore, 
$ B=\begin{pmatrix}
a & 0\\
0 & c
\end{pmatrix}$, where
$ a,c \in R $; the other conditons lead to 
\[
N(A_1)= \lbrace B=\begin{pmatrix}
a & 0\\
0 & c
\end{pmatrix}~|~ a, c \in R, ~~ B \neq aI, ~~ B\neq A_1 \rbrace .
\]
Obviously, $ d(A_1)=|N(A_{1})|=|R|^2 -|R| -1. $\\ 

The proof of $ (ii) $ is similar to $ (i) $ and suppressed.\\

$ (iii) $ Let 
$A_{3} = \begin{pmatrix}
u & v\\
0 & 0
\end{pmatrix}$, where $ u, v \in U(R) $ are fixed. Take 
$ B = \begin{pmatrix}
a & b\\
0 & c
\end{pmatrix}\in T\setminus Z(T) $
 such that $ B\neq A $ and $ AB=BA $. By an easy computation we arrive at \\
 
\[ 
 N(A_3)= \lbrace B=\begin{pmatrix}
a & b\\
0 & a-u^{-1}vb
\end{pmatrix}~|~ a,b \in R, ~~ B \neq aI, ~~ B\neq A_3 \rbrace . 
\]
Since $ u, v $ are fixed, we have again $ d(A_3)=|R|^2 -|R| -1. $\\

$ (iv) $ Let 
$ A_{4} = \begin{pmatrix}
u & n\\
0 & 0
\end{pmatrix} $, with $ u \in U(R) $ and $ n\in X $. Then in order to 
$B = \begin{pmatrix}
a & b\\
0 & c
\end{pmatrix}$ $\in T$ make an edge with $ A_{4} $ $ (AB=BA, B\neq aI, B\neq A_{4})$ we find that 
\[
N(A_4)= \lbrace B=\begin{pmatrix}
a & b\\
0 & c
\end{pmatrix}~|~ a,c \in R,~~ b \in Z(n), ~~ B \neq aI, ~~ B\neq A_4 \rbrace.
\]
Therefore, 
$ d(A_4)=|Z(n)|~|R|^2 -|R| -1. $\\

The proof of $(v)$ is similar to $(iv)$ and omitted.\\

$ (vi) $ Now let 
$A_{6} = \begin{pmatrix}
n_{1} & n_{2}\\
0 & 0
\end{pmatrix}\in $ $T$ with $n_{i} \in X$. Assume that 
$B = \begin{pmatrix}
a & b\\
0 & c
\end{pmatrix}$
is adjacent to $A_{6}$. Then we find that 
$$ N(A_6)= \lbrace  B=\begin{pmatrix}
a & b\\
0 & c
\end{pmatrix}~|~ a,b,c \in R, bn_1+(c-a)n_2=0, B \neq aI, B\neq A_6 \rbrace. $$
The condition $ bn_{1}+(c-a)n_{1}=0 $ means that $ b $ and $ c-a $ satisfy the equation $ n_{1}x+n_{2}y=0 $, and this implies that two of the variables $ a, b, c $ in matrix $B$ are arbitrary and the third one is determined by the given condition. Considering also the cases $B\neq aI, B\neq A_{6}$, if $S$ is the number of all solutions of $ n_{1}x+n_{2}y=0 $ in $R$, we find that
$ d(A_6)=S~|R|^2 -|R| -1, $ as desired.\\

$ (vii) $ Let 
$A_{7} = \begin{pmatrix}
u & n\\
0 & 0
\end{pmatrix}$ with $u \in U(R), n \in X$, and let
$B = \begin{pmatrix}
a & b\\
0 & c
\end{pmatrix}$
be adjacent to $A_{7}$. 
Repeating similar arguments as in previous parts, we observe that \\
$$ N(A_7)= \lbrace B=\begin{pmatrix}
a & u^{-1}(a-c)n\\
0 & c
\end{pmatrix}~|~ a,c \in R, ~~ B \neq aI, ~~ B\neq A_7 \rbrace .$$ 
Therefore,
$ d(A_7)=|R|^2 -|R| -1. $\\

$(viii)$ Proofs of $ (viii) $ and $ (vii) $ are similar, hence omitted.
\qed\\

The first corollary of the theorem above determines the number of edges of the commuting graph of $ 2 \times 2 $ upper triangular matrix rings over a finite field.
\begin{cor}\label{cor2.3}
Let $ F $ be a finite field of order $ p^n $, where $ p $ is a prime and $ n\geqslant 1 $, and let $ G $ be the commuting graph $ \Gamma (Tr (F)) $. Then
$$ |E (G)| =  \dfrac{1}{2} p^n (p^{2n}-1) (p^{2n}-p^n-1).$$
\end{cor}
\proof
Since $ X= \lbrace r\in R\, |\, r\neq 0\,\, \text{and}\,\, r \,\, \text{is a zero divisor} \rbrace =\emptyset $, we infer that $ A_{i} ^{*} =0$ for all $ i\geqslant 4 $. Therefore, by Lemmas \ref{10} , \ref{(2.1)}, Theorem \ref{thm2.2} and (\ref{8}), (\ref{9}) 
we have 
\begin{align*}
|E(G)| & =\dfrac{1}{2} \sum_{i=1} ^{3} d(A_{i} ^{*})\\
&=\dfrac{1}{2} (2p^{n} (p^{n}-1)(p^{2n} -p^{n}-1)+p^{n}(p^{n}-1)^{2}(p^{2n}-p^{n}-1))\\
&=\dfrac{1}{2} p^{n}(p^{2n}-1)(p^{2n}-p^{n}-1).
 \end{align*}
\qed\\

We can also obtain the number of edges of the commuting graph $ \Gamma (Tr (F)) $ when $ R $ is the direct product of two finite fields. The long computation for case of finite direct products of finite fields can be handled similary.
\begin{cor}\label{cor2.4}
Let $ F_1 $, $ F_2 $ be finite fields of order $ p^m $, $ q^n $, respectively, where $ p,q $ are not necessarily distinct primes and $ m,n\geqslant 1 $. Let $ G=\Gamma (Tr (F_1 \times F_2)) $. Put $ R=F_1 \times F_2 $. Then
\begin{align*}
|E (G)| =  \dfrac{|R|}{2}(Y |U| (4+|U| +2 |U_1| + |U_2| )+ Z |U_1| (2+|U_1| )+ W  |U_2| (2+ |U_2|)),
\end{align*}
 where $ |R|=|F_1|~|F_2| = p^m q^n $, $ |U_1|=|U(F_1)|=p^m-1 $,
 $ |U_2| =|U(F_2)|=q^n -1 $, $ |U|=|U(R)|=(p^m-1)(q^n -1) $,  $ Y=p^{2m} q^{2n} -p^m q^n -1 $,
  $ Z=p^{2m} q^{3n} -p^m q^n -1  $,
 and $ W=p^{3m} q^{2n} -p^m q^n -1  $.
\end{cor}
\proof
Let $ |F_{1}|=p^{m} $, $ |F_{2}|=q^{n} $, $ G=\Gamma (Tr(F_{1} \times F_{2})) $,\\ 
$ A_{1}=
\begin{pmatrix}
(u,v) & (0,0)\\
(0,0) & (0,0)
\end{pmatrix} $,
$ A_{2}=
\begin{pmatrix}
(0,0) & (u,v)\\
(0,0) & (0,0)
\end{pmatrix} $,
$ A_{3}=
\begin{pmatrix}
(u_{1},v_{1}) & (u_{2},v_{2})\\
(0,0) & (0,0)
\end{pmatrix} $,\\
where $ u,u_{i} \in U(F_{1}) $ and $ v,v_{i}\in U(F_{2}) $.  
Then, by Lemma \ref{lem2.1} and Theorem \ref{thm2.2}\\
 we have 
\[
d(A_{2} ^{*})=d(A_{1} ^{*})=|A_{1}^{*}||N(A_{1})|=|U||R|(|R| ^{2} -|R|-1);
\]
\[
d(A_{3} ^{*})=|A_{3}^{*}||N(A_{3})|=|U|^{2}|R|(|R| ^{2} -|R|-1).
\]
We have $ A_{4} = \begin{pmatrix}
(x,y) & (0,0) \\
(0,0) & (0,0)
\end{pmatrix}$,
 where $(x,y) \in F_{1} \times F_{2} -U(F_{1})\times U(F_{2})- \lbrace(0,0)\rbrace.$ Therefore, $ A_{4} $ has one of the forms 
$ A_{4,1} =\begin{pmatrix}
(u,0) & (0,0) \\
(0,0) & (0,0)
\end{pmatrix}$, or
$ A_{4,2} = \begin{pmatrix}
(0,v) & (0,0) \\
(0,0) & (0,0)
\end{pmatrix}$,
where $ u \in U(F_{1})=F_{1} \setminus \lbrace 0 \rbrace $ and $ v \in U(F_{2})=F_{2} \setminus\lbrace 0 \rbrace $. By Theorem \ref{thm2.2} ($ iv $) we have 
 \begin{align*}
d(A_{4,1} ^{*})=|A_{4,1}^{*}|d(A_{4,1})& =|U_{1}||R|(|Z(u,0)||R| ^{2} -|R|-1)\\
&=|U_{1}||R| (|F_{2}||R|^{2}-|R|-1);
 \end{align*}
\begin{align*}
d(A_{4,2} ^{*})=|A_{4,2}^{*}|d(A_{4,2})& =|U_{2}||R|(|Z(0,v)||R| ^{2} -|R|-1)\\
&=|U_{2}||R| (|F_{1}||R|^{2}-|R|-1),
 \end{align*}
so that
 \begin{align*}
d(A_{5} ^{*})=d(A_{4}^{*})&=d(A_{4,1} ^{*})+ d(A_{4,2} ^{*})\\
& =|R|(|U_{1}| (|F_{2}||R|^{2}-|R|-1)+|U_{2}| (|F_{1}||R|^{2}-|R|-1)).
 \end{align*}
Now, by definition, $ A_{6} $ has one of the following forms:\\
$$ A_{6,1} = \begin{pmatrix}
(u_{1},0) & (u_{2},0) \\
(0,0) & (0,0)
\end{pmatrix},
 A_{6,2} = \begin{pmatrix}
(u_{1},0) & (0,v_{1}) \\
(0,0) & (0,0)
\end{pmatrix}, $$ 
$$ A_{6,3} = \begin{pmatrix}
(0,v_{1}) & (u_{1},0) \\
(0,0) & (0,0)
\end{pmatrix},
A_{6,4} = \begin{pmatrix}
(0,v_{1}) & (0,v_{2}) \\
(0,0) & (0,0)
\end{pmatrix}, $$ 
where $ u_{i} \in U(F_{1}) $ and $ v_{i} \in U(F_{2}) $. By Theorem \ref{thm2.2} ($ iv $)
\begin{align*}
N(A_{6,1}) =& \lbrace B= \begin{pmatrix}
(a_{1},a_{2}) & (b_{1},b_{2})\\
(0,0) & (c_{1},c_{2})
\end{pmatrix}
| (a_{1},a_{2}), (b_{1},b_{2}),(c_{1},c_{2}) \in R, \\
&\,\,\, (b_{1},b_{2})(u_{1},0)+(c_{1}-a_{1},c_{2}-a_{2})(u_{2},0)=0, B\neq (a_{1},a_{2})I, B\neq A_{6,1}\rbrace .
 \end{align*}
So that $ b_{1} u_{1}+(c_{1}-a_{1})u_{2}=0 $, and $ b_{1}=-u_{1} ^{-1} u_{2} (c_{1}-a_{1}). $
Therefore, $ B $ is of the form 
$\begin{pmatrix}
(a_{1},a_{2}) & (-u_{1} ^{-1} u_{2} (c_{1}-a_{1}), b_{2}) \\
(0,0) & (c_{1},c_{2})
\end{pmatrix}$,
where $ (a_{1},a_{2}), (c_{1},c_{2}) \in R $ and $ b_{2} \in F_{2} $ are all arbitrary. Considering the other restrictions on $ B $, yields 
\[
d(A_{6,1})=|F_{2}||R|^{2}-|R|-1=:Z.
\]
A similar argument shows that 
\[
d(A_{6,4})=|F_{1}||R|^{2}-|R|-1=:W.
\]
Now, consider 
$A_{6,2}= \begin{pmatrix}
(u_{1},0) & (0,v_{1})\\
(0,0) & (0,0)
\end{pmatrix}$.
 By Theorem \ref{thm2.2} ($ iv $) again,
\begin{align*}
N(A_{6,2}) =& \lbrace B= \begin{pmatrix}
(a_{1},a_{2}) & (b_{1},b_{2})\\
(0,0) & (c_{1},c_{2})
\end{pmatrix}
| (a_{1},a_{2}), (b_{1},b_{2}),(c_{1},c_{2}) \in R, \\
&\,\,\, (b_{1},b_{2})(u_{1},0)+(c_{1}-a_{1},c_{2}-a_{2})(0,v_{1})=0\rbrace,
 \end{align*}
so that $ b_{1} u_{1}=0=(c_{2}-a_{2})v_{1} $, whence $ b_{1}=0 $ and $ c_{2}=a_{2} $. These restrictions imply that
$B=\begin{pmatrix}
(a_{1},a_{2}) & (0,b_{2})\\
(0,0) & (c_{1},a_{2})
\end{pmatrix}.$
Considering the other two restrictions on $ B $, we get $d(A_{6,2})=|R|^{2} -|R|-1=:Y.$
A similar computation shows that $d(A_{6,3})=d(A_{6,2})=Y. $ Consequently, by Lemma \ref{lem2.1} we have
\begin{align*}
d(A_{6} ^{*}) = \sum _{i=1} ^{4} d(A_{6,i} ^{*} )& = \sum_{i=1} ^{4} |A_{6,i} ^{*}|d(A_{6,i})\\
& = |U_{1}|^{2} |R| Z + |U_{2}|^{2} |R| W+2|U||R|Y.
\end{align*}
Finally, to find $d(A_{7} ^{*})(=d(A_{8} ^{*})),$ note that $ A_{7} $ can have one of the following forms\\
$A_{7,1}= \begin{pmatrix}
(u,v) & (u_{1},0)\\
(0,0) & (0,0)
\end{pmatrix} $
or
$A_{7,2}= \begin{pmatrix}
(u,v) & (0,v_{1})\\
(0,0) & (0,0)
\end{pmatrix} $,
 where $ u, u_{1} \in U(F_{1}) $ and $ v, v_{1} \in U(F_{2}) $. By Theorem \ref{thm2.2} ($ iiv $) we have $ d(A_{7,1})=d(A_{7,2})= |R|^{2} -|R|-1 =Y. $ Hence, 
\[
 d(A_{7} ^{*})=d(A_{7,1}^{*})+ d(A_{7,2}^{*})=(|U_{1}|+|U_{2}|) |U||R|Y=d(A_{8} ^{*}).
\]
Finally, using Lemma \ref{lem2.1} and the above computation, we find that 
\begin{align*}
|E (G)| = \dfrac{|R|}{2}(Y |U| (4+|U| +2 (|U_1| + |U_2| )+ Z |U_1| (2+|U_1| )+ W |U_2| (2+ |U_2|)).
\end{align*}
\qed
\begin{prop}\label{pro2.5}
Let $ R $ and $ S $ be isomorphic noncommutative finite rings with identity.\\
Then $ \Gamma (R) \cong \Gamma (S) $.
\end{prop}
\proof
Let $ \varphi: R \longrightarrow S $ be a ring isomorphism, and let $ G_{1} = \Gamma (R) $ and $ G_{2}=\Gamma (S)$. Then $ \overline{\varphi} : G_{1} \longrightarrow G_{2}$ given by $ \overline{\varphi}(r)=\varphi (r) $ is easily seen to be a graph isomorphism.
\qed
\begin{ex}\label{exa2.6}
$(a)$  Let $ G_1 = \Gamma (Tr (\mathbb{Z}_{10})) $. Then, by the proposition above,
$$ G_1\cong  \Gamma (Tr(\mathbb{Z}_{2} \times \mathbb{Z}_{5})). $$
Now, by Corollary \ref{cor2.4} we have $ |R|=10 $, $ Y=89 $, $ Z=489 $, $W=189 $, $ |U|=4 $, $ |U_1|=1 $, and $  |U_2|=4 $, so that $| E(G_1)| =62055 $.\\

$(b)$ Let $ G_2 = \Gamma (Tr(\mathbb{Z}_{5}\times \mathbb{Z}_{5})) $. Then, by Corollary \ref{cor2.4}   we  have\\
$  |R|=25 $, $ Y=599 $, $ Z=W=3099 $, $ |U|=16 $, and $ |U_1|=|U_2|=4 $. So that $ |E(G_2)|=6172000. $
\end{ex}
\begin{rem} \label{R2.8}
For the case when $R$ and $S$ are finite fields, the converse of Proposition \ref{pro2.5} holds: For if $ \Gamma (Tr(R)) \cong \Gamma (Tr(S)) $, then $ |R|^{3}-|R|=|E(\Gamma (|R|)|=|E(\Gamma |S|)| = |S|^{3} -|S|$. Since $ |R| $ and $ |S| $ are power of primes   , we infer that $ |R|=|S|=p^{n} $ for some prime $ p $ and $ n\geqslant 1 $. Hence, by an standard theorem in finite field theory, $ R $ and $ S $ are splitting feilds of the polynomial $ x^{p^{n}}-x \in \mathbb{Z}_{p} [x] $, and as such $ R\cong S $.
\end{rem}
\indent However, we do not know if the converse of Proposition \ref{pro2.5} is true in general.
\begin{cor}\label{cor2.7}
Let $ T= Tr(\mathbb{Z} _n)$, $ n\geqslant2 $, and let $ G=\Gamma (T) $. For $ A _i ,A _i ^* $ given in (\ref{6}), (\ref{8}) and the notation defined in Lemma \ref{lem2.1}, 
we have
\begin{align*}
d(A_7 ^*)=d(A_8 ^*)= n \varphi (n) (n -\varphi (n) -1) \varphi (n^2 -n-1).
\end{align*}
\end{cor}
\proof
Let $A_{7} = \begin{pmatrix}
u & m\\
0 & 0
\end{pmatrix}$, where $u\in U$ and $m\in X$. Noting that $|U(\mathbb{Z}_{n})|=\varphi(n) $ and $ |X|=n-\varphi(n)-1 $, by Lemma \ref{lem2.1} and Theorem \ref{thm2.2} ($ vii $) we have
$$d(A_{7} ^{*})=|A_{7} ^{*}|d(A_{7})=n\varphi (n)|X|(n^{2}-n-1)=n\varphi (n)(n-\varphi(n)-1)(n^{2}-n-1).$$
Clearly, $ d(A_{8}^{*})=d(A_{7}^{*}) $.
\qed\\

To avoid tedious computation, in the next corollary which computes the remained $d(A_{i}) ^{*^{,}}$s,  the ring $ R $ is more restricted.
\begin{cor}\label{cor2.8}
Let $ T= Tr(\mathbb{Z} _{p^n})$, where $ p $ is a prime and $ n \geqslant 2 $. Put $ G=\Gamma (T) $. Then
\begin{enumerate}
\item[$(i)$ ] $ d(A_4 ^*)=d(A_5 ^*)=(p-1) p^{2n-1} (\sum _ {r=1} ^{n-1} (p^{2n} -p^{n-r} -p^{-r})); $\\
\item[$(ii)$] $ d(A_6 ^*)= 2 \alpha - \beta , ~~\text{where} $\\

\hspace*{0.5cm}$ \alpha = (p-1) p^{2n-1} (\sum _ {r=1} ^{n-1}((p^{3n} +p^{n}+1) p^{-r} -(p^{2n}+p^n) p^{-2r} - p^{2n})), $\\
\hspace*{.9cm}$ \beta = (p-1)^{2} \sum _ {r=1} ^{n-1} ( p^{3n-2r-2 } ( p^{2n+1} -p^n -1)). $
 
\end{enumerate}
\hspace*{0.5cm} Moreover, the clique number $w(G)$ of $G$ is $p^{n}(p^{n}-1)$ which happens twice; and the independence number $ \alpha(G) $ of $ G $ is $ p^{n}(p^{2n}-p^{2n-2}) $.
\end{cor}
\proof
($ i $) We have $ A_{4} = \begin{pmatrix}
m & 0\\
0 & 0
\end{pmatrix} $, for some $ m \in $ $ (p) \setminus \lbrace 0\rbrace $. Let $ m=p^{r} \alpha $ for some fixed $ 1\leqslant r \leqslant n-1 $ and $ \alpha = \alpha (p) = \alpha _{0}+ \alpha_{1} p +...+\alpha_{n-r-1} p^{n-r-1} \in U(\mathbb{Z} _{p^{n}} )  $, which is polynomial of degree at most $ n-r-1 $ on $ p $ such that $ 0 \leqslant \alpha_{i} \leqslant p-1 $, except that $ \alpha_{0} \neq 0$. Put $ A_{4r} = \begin{pmatrix}
p^{r} \alpha & 0\\
0 & 0
\end{pmatrix}$,
and assume that $ B = \begin{pmatrix}
a & b\\
0 & c
\end{pmatrix}$$ \in Tr(\mathbb{Z} _{p^{n}})$ satisfies $ AB=BA. $ This amounts to $ p^{r} \alpha b=0; $ hence $ p^{r} b=0 $, and $ b \in Z(p^{r})=(p^{n-r}). $ Therefore, 
\[
N(A_{4r})= \lbrace B=\begin{pmatrix}
a & b\\
0 & c
\end{pmatrix} | a,c \in \mathbb{Z} _{p^{n}} , b \in (p^{n-r}), B \neq aI, B\neq A_{4r} \rbrace ;
\]
\[
d(A_{4r})=|(p^{n-r})| p^{2n} - p^{n} -1=p^{2n+r} - p^{n} -1.
\]
Now, by Lemma \ref{lem2.1} we have
\begin{align} \label{3.1}
d(A_{4r} ^{*})= |A_{4r} ^{*}| d(A_{4r})=(p-1)p^{n-r-1} p^{n} (p^{2n+r} -p^{n}-1),
\end{align}
in which $ (p-1)p^{n-r-1} $ is the number of choices of $ \alpha . $ Noting that for $ r \neq s $, $ A_{4r} ^{*} \bigcap A_{4s} ^{*} = \emptyset $, adding (\ref{3.1}) over all $ 1\leqslant r \leqslant n-1, $ we get 
\begin{align*}
d(A_{4} ^{*}) =\sum _{r=1} ^{n-1} d(A_{4r} ^{*})& =\sum _{r=1} ^{n-1} (p-1) p^{2n-r-1} (p^{2n+r} -p^{n}-1)\\
& =  (p-1) p^{2n-1} \sum _{r=1} ^{n-1}  (p^{2n} -p^{n-r}-p^{-r})= d(A_{5}).
\end{align*}

($ ii $) Put $ A_{6,r,s}= \begin{pmatrix}
p^{r} x & p^{s} y\\
0 & 0
\end{pmatrix} $, where $ 1 \leqslant r,s \leqslant n-1 $ and $ x , y \in U(\mathbb{Z} _{p^{n}} ) $ with $ \deg x \leqslant n-r-1 $ and $ \deg y \leqslant n-s-1 $ are all fixed. Assume that $ r \leqslant s $, and $ B= \begin{pmatrix}
a & b\\
0 & c
\end{pmatrix}$ $\in Tr(\mathbb{Z} _{p^{n}})$ satisfies $AB=BA$. Therefore, $ p^{r} x b + p^{s} y c -p^{s} x y =0 $, so that $ p^{r}(b+x ^{-1} y p^{s-r} (c-a))=0 $, hence $ b+x^{-1} y p^{s-r}(c-a) \in (p^{n-r})$. Note also that the latter relation implies that $AB=BA$. Consequently, 
\begin{align}\label{3.2}
 N(A_{6,r,s})= \lbrace B=\begin{pmatrix}
a & b\\
0 & c
\end{pmatrix} | a,c \in R, b+x ^{-1} y p^{s-r}  (c-a) \in (p^{n-r} ) , B \neq aI, B\neq A_{6,r,s} \rbrace ;
\end{align}
$ d(A_{6,r,s})=p^{2n+r}-p^{n}-1 $, for fixed $ r, s,x ,y $, and $ r\leqslant s $. Adding (\ref{3.2}) over all $ 1 \leqslant r \leqslant s \leqslant n-1 $ and possible $ x , y, $ we arrive at
\begin{align*}
d(A_{6, r\leqslant s}) & = \sum _{r=1} ^{n-1} \sum _{s=r} ^{n-1} (p-1)^{2} p^{n} p^{n-r-1}  p^{n-s-1} (p^{2n+r}- p^{n}-1) \\
&=(p-1)^{2} \sum _{r=1} ^{n-1} p^{2n-r-1} (p^{2n+r}-p ^{n} -1) \sum _{ s=r} ^{n-1} p^{n-s-1} \\
& = (p-1) p^{2n-1} ( \sum _{r=1} ^{n-1} (p^{3n}+p^{n}+1)p^{-r} -(p^{2n}+p^{n})(p^{-2r}-p^{2n}))=: \alpha.
\end{align*}
 A simple observation shows that $ d(A_{^{6,r,s(s\leqslant r)}}  ^{*})=d(A_{^{6,r,s(r\leqslant s)}}  ^{*})=\alpha $. Now let $\beta =d(A_{6,r,r} ^{*})  $, where  $A_{6,r,r}=\begin{pmatrix}
p^{r} x & p^{r} y\\
0 & 0
\end{pmatrix}$. Then putting $ r=s $ in $ \alpha $, we get
\[
\beta = d(A_{6,r,r} ^{*})=(p-1)^{2} \sum _{r=1} ^{n-1} p^{3n-2r-2} (p^{2n+r}-p^{n}-1).
\]
Therefore, $ d(A_{6} ^{*})=d(A_{6,r\leqslant s} ^{*}) +d(A_{6, s\leqslant r} ^{*})-d(A_{6,r=s} ^{*})=2\alpha -\beta .$\\
The clique number of $ G $ is obtained in a more general condition in Theorem \ref{thm2.10}, in which $ G=\Gamma (Tr( \mathbb{Z} _{m})) $, where $ m $ is not a prime. To prove the last part, consider the following subsets of the set $ S= \lbrace A_{i} ^{*} | 1\leqslant i \leqslant 8 \rbrace $:
\[
S_{1}= \lbrace A_{1} ^{*},A_{2} ^{*},A_{3} ^{*},A_{7} ^{*},A_{8} ^{*} \rbrace,
S_{2}= \lbrace A_{1} ^{*},A_{2} ^{*},A_{6} ^{*} \rbrace,
S_{3}= \lbrace A_{1} ^{*},A_{3} ^{*},A_{5} ^{*} \rbrace,
S_{4}= \lbrace A_{5} ^{*},A_{4} ^{*} \rbrace.
\]
It is not hard to see that for each $ 1 \leqslant j \leqslant 4 $, the elements in $ S_{j} $ are mutually independent; that is, no pair of matrices from distinct elements of $ S_{j} $ are adjacent. Moreover, $ S_{j}, 1\leqslant j \leqslant 4 $, are the only subsets of $ S $ with this property. We claim that $ \sum _{A_{i} ^{*} \in S_{1}} |A_{i} ^{*} | $
is maximal. To see this, first note that:\\
$ |A_{1} ^{*}|=p^{n}(p^{n}-p^{n-1})=|A_{2}^{*}| $, $  |A_{3} ^{*}|=p^{n}(p^{n}-p^{n-1})^{2} $, $|A_{4} ^{*}|=|A_{5} ^{*}|=p^{n}|X| $, $ |A_{6} ^{*}|=p^{n}|X|^{2} $, and $ |A_{7} ^{*}|=|A_{8} ^{*}|=p^{n}(p^{n}-p^{n-1})|X| $, where $ |X|=p^{n-1} -1 $.\\
We show that $ \sum_{A_{i} ^{*} \in S_{1}} |A_{i} ^{*}| \geqslant \sum_{A_{i} ^{*} \in S_{1}} |A_{i} ^{*}| $. To do this, it suffices to show that $ |A_{3} ^{*}|+2 |A_{7} ^{*}| > |S_{6} ^{*}|$. This is equivalent to $ ((p^{n}-p^{n-1})+|X|^{2}) > 2|X|^{2} $, or $ \dfrac{(p^{n}-1)^{2}}{(p^{n-1}-1)^{2}} > 2 $.
Since $ p\geqslant 2 $, we have $ p^{n}-1 \geqslant 2p^{n-1} -1 > 2p^{n-1} -2 $, so that $ \dfrac{p^{n}-1}{p^{n-1}-1} > 2 $, hence $ \dfrac{(p^{n}-1)^{2}}
{(p^{n-1}-1)^{2}} > 2 $, as desired. Clearly, 
$ \sum_{i \in S_{1}} |A_{i} ^{*}| > \sum_{i \in S_{3}} |A_{i} ^{*}|  $ and $ \sum_{i \in S_{1}} |A_{i} ^{*}| > \sum_{i \in S_{4}} |A_{i} ^{*}|  $. Therefore, the independence number of $ G $ is $ \alpha (G)= \sum_{i \in S_{1}} |A_{i} ^{*}| =2 |A_{1} ^{*}|+|A_{3} ^{*}|+2 |A_{7} ^{*}|=p^{n}(p^{2n}-p^{2n-2}) $, as desired.
\qed\\

Now, using Theorem \ref{thm2.2}, Corollaries \ref{cor2.3}, \ref{cor2.7}, \ref{cor2.8}, and Handshaking Lemma \ref{10}, one can easily find the number of edges of the graph $ G=\Gamma(Tr (\mathbb{Z} _{p^n})) $, where $ p $ is a prime and $ n\geqslant2 $:
\begin{align}\label{m2.2}
|E(G)| = \dfrac{1}{2} \sum _{i=1} ^{8} d(A_i ^*).
\end{align}
\begin{ex}\label{exa2.8}
 $(a)$ Let $ p $ be a prime, and let $ G_1 = \Gamma (Tr (\mathbb{Z} _{p^2})) $. Then 
\begin{align*}
d(A_1 ^*)&=d(A_2 ^*) = p^3 (p-1) (p^4 -p^2 -1),~ d(A_3 ^*)=p^4 (p-1) ^2 (p^4 -p^2 -1),\\
d(A_4 ^*)&=d(A_5 ^*) = p^2 (p-1) (p^5 -p^2 -1),~ d(A_6 ^*)=p^2 (p-1) ^2 (p^5 -p^2 -1),\\
d(A_7 ^*)&=d(A_8 ^*) =  p^3 (p-1)^2 (p^4 -p^2 -1).
\end{align*}Therefore,
 \begin{align*}
 |E(G_1)| = \dfrac{1}{2} \sum _{i=1} ^ 8 d(A_i ^*)=\dfrac{1}{2}p^2 (p^2-1) (p^6 + p^5 -p^4  -2p^2 -1).
 \end{align*}
\indent $(b)$ Let $ G_2 =\Gamma(Tr (\mathbb{Z}_8)) $. Then
 \begin{align*}
 d(A_1 ^*)=&d(A_2 ^*) =1760,~~ d(A_3 ^*)=7040,~~d(A_4 ^*)=d(A_5 ^*)=3880,\\
& d(A_6 ^*)=9592,~~ d(A_7 ^*)=d(A_8 ^*)=5280.
 \end{align*}
 Hence, $  |E(G_2)| = \dfrac{1}{2} \sum _{i=1} ^ 8 d(A_i ^*)=19236 $.
 \end{ex}
\indent Now, we are ready to state more results about the commuting graph of $ 2 \times 2 $ upper triangular matrix rings over finite fields.\\
 
First, we have a definition: let $ F $ be a finite field of order $ p^n (n\geqslant 1) $. Fix $ x \in U(F) $, define
 $ A_{3x} = \begin{pmatrix}
u & xu\\
0 & 0
\end{pmatrix} $, where $ u $ is also any fixed element
 in $ U(F) $, and put
\[
 A_{3x}^* = \lbrace  \begin{pmatrix}
v+z & xv\\
0 & z
\end{pmatrix} ~|~ v \in U(F), ~~ z \in F \rbrace .
\]
\indent  Note that
 \begin{align*}
 |A^* _1| =|A^* _2|=|A^* _{3x}|=p^n (p^n -1) ~~~(\text{for~all}~~ x \in U(F)),
 \end{align*}
 and $ \sum _{x \in U(F)} |A^* _{3x}| =|A^* _{3}| =p^n (p^n -1)^2 $.\\

Let $ F $ be a finite ring. In 2004, Akbari et al. \cite{AGHM} have shown that the number of connected components (cliques) of $ \Gamma (M_{2} (F)) $ is $ |F|^{2}+|F|+1 $, and each of them has $ |F|^{2}-|F| $ vertices. In the following theorem we show that $ \Gamma (T_{2} (F)) $ has $ |F|^{2}+1 $ cliques each with $ |F|^{2} -|F| $ vertices.
\begin{thm}\label{thm2.9}
Let $ F $ be a finite field of order $ p^n (n\geqslant1) $, $ T=Tr(F) $, and $ G=\Gamma (T) $. Then\\
\indent $(i)$ Each of the following mutually disjoint $ p^n +1 $ subsets $ A^* _1 $, $ A^* _2 $, and $ A^* _{3x} $ (when $ x $ ranges over $ U(F) $) of $ G $ induces a connected component (clique) of $ G $. These have the same number of vertices $ p^{2n} -p^n $, the same number of edges 
$ \begin{pmatrix}
p^{2n}-p^{n}\\
2
\end{pmatrix} $,
 and the same diameter $ 1 $.\\
\indent $(ii)$ The clique number $ w(G) $ of $ $ G $ $ is $ p^{2n} -p^n $; $ \delta (G) =\Delta (G) = p^{2n} -p^n -1  $; $ \diam G=\infty $; $ G $ is $ (p^{2n} -p^n -1) $-regular, and the independence number $ \alpha (G)  $ of $ G $ is $ p^n+1 $.
\end{thm}
\proof
($ i $) Obviously, $ V(G) $ is the disjoint union of the following $ |F|+1 $ subsets 
\begin{align} \label{3.3}
A_{1}^{*}, A_{2} ^{*}, A_{3x} ^{*} (x \in U(F)),
\end{align}
with the same cardinality $ p^{2n} - p^{n} $. Let $ A $ be any of these subsets. Then, by Theorem \ref{thm2.2} ($ i $)-($ iii $)
we have $ A^{*}=N(A) \cup \lbrace  A \rbrace$; hence for all $ \beta_{i}, \beta_{j} \in A_{i} ^{*} $, $ \beta_{i} \beta_{j} = \beta_{j} \beta_{i}, $ concluding that subgraph $ G[A^{*}]=(A^{*}, \lbrace  (\beta_{i}, \beta_{j})| \beta_{i}, \beta_{j} \in A^{*} \rbrace )$ is a (complete) connected subgraph with 
$ \begin{pmatrix}
p^{2n}-p^{n}\\
2
\end{pmatrix} $
edges whose diameter is 1.\\

($ ii $) Since the subsets in (\ref{3.3}) make a partition of $ V(G), $ we see that $ w(G)= p^{2n}-p^{n} $, and $ \delta(G)=\Delta(G) =d(A)= p^{2n}-p^{n}-1 $, where $ A $ is any of the subsets in (\ref{3.3}). By definition, no pairs $ A, A^{'} $ from two different subsets in (\ref{3.3}) are adjacent; meaning that $ \diam G =\infty  $, each $ G[A^{*}] $ is a component, and that the independence number $ \alpha (G) $ of $ G $ is the number of subsets in (\ref{3.3}), which is $ |F|+1 $. Clearly, $ G $ is $ (p^{2n}-p^{n}-1) $-regular. This completes the proof.
\qed\\

Our final result may be of some interests:
\begin{thm}\label{thm2.10}
Let $ G=\Gamma(Tr (\mathbb{Z}_m)) $, where $ \mathbb{Z}_m $ is not a field. Then $ G $ is connected, $ \diam G=3 $, and the clique number of $ G $ is $ m(m-1) $ which happens twice.
\end{thm}
\proof
Let $ T=Tr(\mathbb{Z} _{m}), U=U(\mathbb{Z} _{m}), X=\mathbb{Z}_{m} -U(\mathbb{Z}_{m})-\lbrace 0\rbrace, $ and 
$ G=\Gamma (T) $. Using Remark \ref{1.1} and the fact that $ V(G)= \bigcup  ^{8}  _{i=1} A_{i} ^{*} $, it suffices to consider the following mutually disjoint sets $ \overline{A_{i}} $ of matrices: \\

\noindent $ \overline{A_{1}}= \lbrace \begin{pmatrix}
u & 0 \\
0 & 0
\end{pmatrix} | u \in U \rbrace $, $\overline{A_{2}}= \lbrace \begin{pmatrix}
0 & u \\
0 & 0
\end{pmatrix} | u \in U \rbrace$, $\overline{A_{3}}= \lbrace \begin{pmatrix}
u & v \\
0 & 0
\end{pmatrix} | u,v \in U \rbrace $,\\ $ \overline{A_{4}}= \lbrace \begin{pmatrix}
n & 0 \\
0 & 0
\end{pmatrix} | n \in X \rbrace $, and so on. \\

\indent It is straightforward to verify that:
\begin{itemize}
\item[$ (i) $]
 For all $A,B \in \overline{A_{1}} \bigcup \overline{A_{4}} $  (resp. $A,B \in \overline{A_{2}} \bigcup \overline{A_{5}} ),$ the pair $A, B$ are adjacent.\\
\item[$ (ii) $]
The sets in $ \lbrace \overline{A_{1}} , \overline{A_{2}} ,\overline{A_{6}} \rbrace $ are mutually independent (i.e, no pairs from different sets are adjacent). Moreover, the same conclusion holds for the following sets: 
$$ \lbrace \overline{A_{1}} , \overline{A_{2}} ,\overline{A_{3}}, \overline{A_{7}}, \overline{A_{8}}  \rbrace, \lbrace \overline{A_{1}} , \overline{A_{3}} ,\overline{A_{5}} \rbrace, \lbrace \overline{A_{3}}, \overline{A_{4}}\rbrace. $$
\item[$ (iii) $]
Each $\begin{pmatrix}
n & 0 \\
0 & 0
\end{pmatrix} \in \overline{A_{4}}$ incidents to $\begin{pmatrix}
n_{1} & t \\
0 & 0
\end{pmatrix} \in \overline{A_{6}}$ for all $ n_{1} \in X $ and $ t \in Z(n) $. Conversely, each  $\begin{pmatrix}
n_{1} & n_{2} \\
0 & 0
\end{pmatrix} \in \overline{A_{6}}$ incidents to $\begin{pmatrix}
t & 0 \\
0 & 0
\end{pmatrix} \in \overline{A_{4}}$ for all $ t \in Z(n_{2}) $.\\
\item[$ (v) $]
Each 
 $\begin{pmatrix}
n & 0 \\
0 & 0
\end{pmatrix} \in \overline{A_{4}}$
incidents to 
$\begin{pmatrix}
u & t \\
0 & 0
\end{pmatrix} \in \overline{A_{7}}$
for all $ u \in U,\, t \in Z(n)$;
and conversely, each 
$\begin{pmatrix}
u & n \\
0 & 0
\end{pmatrix} \in \overline{A_{7}}$
incidents to 
$\begin{pmatrix}
t & 0 \\
0 & 0
\end{pmatrix} \in \overline{A_{4}}$
for all $ t \in Z(n) $.
Further, a similar conclusion holds for $ \overline{A_{5}},\, \overline{A_{8}} $, and for $ \overline{A_{4}},\, \overline{A_{8}} $.\\
\item[$ (vi) $]
 Each
$\begin{pmatrix}
n & 0 \\
0 & 0
\end{pmatrix} \in \overline{A_{4}}$ 
(resp. $\begin{pmatrix}
0 & n \\
0 & 0
\end{pmatrix} \in \overline{A_{5}}$) incidents to 
$\begin{pmatrix}
t & 0 \\
0 & 0
\end{pmatrix} \in \overline{A_{5}}$ (resp. 
$\begin{pmatrix}
t & 0 \\
0 & 0
\end{pmatrix} \in \overline{A_{4}}$)
for all $ t \in Z(n)  $.\\

By $ (iii) $, $ \overline{A_{3}} $ is independent with any of $ \overline{A_{1}},\, \overline{A_{2}},\, \overline{A_{4}},\, \overline{A_{5}},\, \overline{A_{7}} $
and $ \overline{A_{8}} $.
Moreover, it is easy to see that the converse of $ (iv) $ is not necessarily true. However, we can construct a path from an arbitrary element of $ \overline{A_{1}} $ to an arbitrary element of $ \overline{A_{3}} $: Take 
$\begin{pmatrix}
u & v \\
0 & 0
\end{pmatrix} \in \overline{A_{3}}$,
and let $ u^{'} \in U $, $ n \in X $ and $ n_{1} \in Z(n) $ be arbitrary. Then, by $ (i) $,
$\begin{pmatrix}
u^{'} & 0 \\
0 & 0
\end{pmatrix}$
incidents to 
$\begin{pmatrix}
n & 0 \\
0 & 0
\end{pmatrix} \in \overline{A_{4}}$; by $ (iii) $,
$\begin{pmatrix}
n & 0 \\
0 & 0
\end{pmatrix} $
incidents to 
$\begin{pmatrix}
n_{1} & u^{-1} v n_{1} \\
0 & 0
\end{pmatrix} \in \overline{A_{6}}$,
and by $ (iv) $,
$\begin{pmatrix}
n_{1} & u^{-1} v n_{1} \\
0 & 0
\end{pmatrix} $ incidents to
$\begin{pmatrix}
u & v \\
0 & 0
\end{pmatrix}$, as desired.\\
\end{itemize}
Now, by  ($ i $),  ($ iii $)-($ vi $)  and the argument above, $ G $ is connected and $ \diam G=3 $. Moreover, one can verify that the subgraph
$G[A_{1} ^{*} \bigcup A_{4} ^{*}]  $ is a clique with (maximal) size
\[
|A_{1} ^{*}|+|A_{4} ^{*}|=m|U|+m|X|=m(m-1)=w(G).
\]
\indent A similar argument holds for 
$ G[A_{2} ^{*} \bigcup A_{5} ^{*}] $.
This completes the proof.
\qed\\
\\

Corollary \ref{cor2.8} and the theorem above suggest the following question:\\
\textbf{Question}:\textit{ What is the independence number of $ \Gamma (Tr(\mathbb{Z}_{n})) $, $ n\geqslant 2 $?}\\

We guess that the answer is $ n \varphi(n)(2n- \varphi (n)). $

\section*{Disclosure statement}
No conflict of interest was reported by the authors.

 %-----------------------------------------------------------

\end{document}